\documentclass[12pt,reqno]{amsart}
\usepackage{amsmath, amssymb, graphicx}
\usepackage{epsfig}
\usepackage[usenames,dvipsnames,svgnames]{xcolor}
\usepackage{tkz-euclide}
\usetikzlibrary{calc,intersections}
\usetkzobj{all}
\definecolor{fondpaille}{cmyk}{0,0,0.1,0}
\renewcommand{\Re}{\operatorname{Re}}
\renewcommand{\Im}{\operatorname{Im}}
\renewcommand{\min}{\operatorname{min}}
\renewcommand{\max}{\operatorname{max}}
\renewcommand{\inf}{\operatorname{Inf}}
\renewcommand{\sup}{\operatorname{Sup}}
\newcommand{\s}{{\sigma}}

\renewcommand{\i}{\infty}
 
\renewcommand{\b}{\beta}
\newcommand{\e}{\epsilon}
\renewcommand{\d}{{\delta}}
\newcommand{\g}{\gamma}
\newcommand{\G}{\Gamma}

\renewcommand{\L}{\Lambda}

\renewcommand{\r}{{\rho}}
\renewcommand{\t}{\theta}

\newcommand{\twopartdef}[4]
{
	\left\{
		\begin{array}{ll}
			#1 & \mbox{if } #2 \\
			#3 & \mbox{if } #4
		\end{array}
	\right.
}
\renewcommand{\(}{\left\(}
\renewcommand{\)}{\right\)}
\renewcommand{\[}{\left\[}
\renewcommand{\]}{\right\]}

\numberwithin{equation}{section}
 \theoremstyle{plain}
\newtheorem{theorem}{Theorem}[section]

   \makeatletter
\def\proof{\@ifnextchar[{\@oproof}{\@nproof}}
\def\@oproof[#1][#2]{\trivlist\item[\hskip\labelsep\textit{#2 Proof of\
#1.}~]\ignorespaces}
\def\@nproof{\trivlist\item[\hskip\labelsep\textit{Proof.}~]\ignorespaces}

\makeatother

\begin{document}
\title[Monotonicity results for Dirichlet L-functions]{Monotonicity results for Dirichlet L-functions}
\author{Atul Dixit}
\address{Department of Mathematics, Tulane University, New Orleans, LA 70118, USA}
\email{adixit@tulane.edu}
\author{Arindam Roy}
\address{Department of Mathematics, University of Illinois, 1409 West Green
Street, Urbana, IL 61801, USA} \email{roy22@illinois.edu}
\author{Alexandru Zaharescu}
\address{Simion Stoilow Institute of Mathematics of the
Romanian Academy, P.O. Box 1-764, RO-014700 Bucharest, Romania and Department of Mathematics, University of Illinois, 1409 West Green
Street, Urbana, IL 61801, USA} \email{zaharesc@illinois.edu}
\begin{abstract}
We present some monotonicity results for Dirichlet $L$-functions associated to real primitive characters. We show in particular that these Dirichlet $L$-functions are far from being logarithmically completely monotonic. Also, we show that, unlike in the case of the Riemann zeta function, the problem of comparing the signs of $\frac{d^k}{ds^k}\log L(s,\chi)$ at any two points $s_1, s_2>1$ is more subtle.
\end{abstract}
\maketitle
\section{Introduction}

A function $f$ is said to be completely monotonic on $[0,\infty)$ if $f\in C[0,\infty)$, $f\in C^{\infty}(0,\infty)$ and $(-1)^{k}f^{(k)}(t)\geq 0$ for $t>0$ and $k=0, 1, 2\cdots$, i.e., the successive derivatives alternate in sign. The following theorem due to S.N.~Bernstein and D.~Widder gives a complete characterization of completely monotonic functions \cite[p.~95]{wawi}:\\

\textit{
A function $f:[0,\infty)\to [0,\infty)$ is completely monotonic if and only if there exists a non-decreasing bounded function $\g$ such that $f(t)=\int_{0}^{\infty}e^{-st}d\g(s)$.
}\\

Lately, the class of completely monotonic functions have been greatly expanded to include several special functions, for example, functions associated to gamma and psi functions by Chen \cite{chen}, Guo, Guo and Qi \cite{qgg} and quotients of $K$-Bessel functions by Ismail \cite{ism}. A conjecture that certain quotients of Jacobi theta functions are completely monotonic was formulated by the first author and Solynin in \cite{dixsoly}, and slightly corrected later by the present authors in \cite{drz1}. Certain other classes of such functions were introduced by Alzer and Berg \cite{alzberg}, Qi and Chen \cite{qichen}. Completely monotonic functions have applications in diverse fields such as probability theory \cite{kimberling}, physics \cite{berg}, potential theory \cite{bergforst}, combinatorics \cite{ball} and numerical and asymptotic analysis \cite{frenzen}, to name a few.

A close companion to the class of completely monotonic functions is the class of logarithmically completely monotonic functions. This was first studied, although implicitly, by Alzer and Berg \cite{alzberg2}. A function $f:(0, \i)\to (0,\i)$ is said to be logarithmically completely monotonic \cite{ber} if it is $C^{\i}$ and $(-1)^k[\log f(x)]^{(k)}\geq 0\mbox{, for }k=0, 1, 2, 3, \cdots$. Moreover, a function is said to be strictly logarithmically completely monotonic if $(-1)^k[\log f(x)]^{(k)}>0$. The following is true:\\

\textit{
Every logarithmic completely monotonic function is completely monotonic.
}\\

The reader is referred to Alzer and Berg \cite{alzberg2}, Qi and Guo \cite{qg}, and  Qi, Guo and Chen \cite{qgc} for proofs of this statement. 

One goal of this paper is to study the Dirichlet $L$-functions from the point of view of logarithmically complete monotonicity. For Re $s>1$, the Riemann zeta function is defined by 
\begin{equation*}
\zeta (s)=\sum_{n=1}^{\i}\frac{1}{n^s}.
\end{equation*} 
Consider $s>1$. Since $\log\zeta(s)>0$ and 
\begin{equation*}
(-1)^k\frac{d^k}{ds^k}\log\zeta(s)=(-1)^k\frac{d^{k-1}}{ds^{k-1}}\left(\frac{\zeta'(s)}{\zeta(s)}\right)=\sum_{n=1}^{\infty}\frac{\Lambda(n)(\log n)^{k-1}}{n^{s}},
\end{equation*}
where $\Lambda(n)\geq 0$ is the von Mangoldt function, $(-1)^k\frac{d^k}{ds^k}\log\zeta(s)>0$ for all $s>1$. This implies that $\zeta(s)$ is a logarithmically completely monotonic function for $s>1$ (in fact, strictly logarithmically completely monotonic). But this approach fails in the case of $L(s,\chi)$ with $s>1$ and $\chi$, a real primitive Dirichlet character modulo $q$, since
\begin{equation*}
(-1)^k\frac{d^k}{ds^k}\log L(s,\chi)=(-1)^k\frac{d^{k-1}}{ds^{k-1}}\left(\frac{L'(s,\chi)}{L(s,\chi)}\right)=\sum_{n=1}^{\infty}\frac{\chi(n)\Lambda(n)(\log n)^{k-1}}{n^{s}}
\end{equation*}
may change sign for different values of $s$ as $\chi(n)$ takes the values $-1, 0$ or $1$. Hence, we need to consider a different approach for studying $L(s,\chi)$ in the context of logarithmically complete monotonicity. This naturally involves studying the zeros of derivatives of $\log L(s,\chi)$.
 
There have been several studies made on the number of zeros of $\zeta^{(k)}(s)$ and $L^{(k)}(s,\chi)$, one of which dates back to Spieser \cite{spe}, who showed that the Riemann Hypothesis is equivalent to the fact that $\zeta'(s)$ has no zeros in $0<$ Re $s<1/2$. Spira \cite{sp1} conjectured that 
\begin{equation*}
N(T)=N_{k}(T)+\left[\frac{T\log 2}{2\pi}\right]\pm 1,
\end{equation*}
where $N_{k}(T)$ denotes the number of zeros of $\zeta^{(k)}(s)$ with positive imaginary parts up to height $T$, and $N(T)=N_{0}(T)$. Berndt \cite{bern} showed that for any $k\geq 1$, as $T\to\infty$,
\begin{equation*}
N_{k}(T)=\frac{T\log T}{2\pi}-\left(\frac{1+\log 4\pi}{2\pi}\right)T+O(\log T).
\end{equation*}
Levinson and Montgomery \cite{lemo} proved a quantitative result implying that most of the zeros of $\zeta^{(k)}(s)$ are clustered about the line Re $s=1/2$ and also showed that the Riemann Hypothesis implies that $\zeta^{(k)}(s)$ has at most finitely many non-real zeros in Re $s<1/2$. Their results were further improved by Conrey and Ghosh \cite{cogh}. Analogues of several of the above-mentioned results for Dirichlet $L$-functions were given by Yildirim \cite{yil2}. Our results in this paper are related to the zeros of $\log L(s,\chi)$ and its derivatives. 

Throughout the paper, we assume that $s$ is a real number and $\chi$ is a real primitive Dirichlet character modulo $q$. Let $F(s,\chi):=\log L(s,\chi)$, and for $s>1$, define 
\begin{equation}
A_{\chi, k}:=\{s: F^{(k)}(s,\chi)=0\}.
\end{equation}
Then we obtain the following result:
\begin{theorem}\label{den}  
Let $\chi$ be a real primitive character modulo $q$ and $L(s,\chi)\neq 0$ for $0<s<1$. Then there exists a constant $c_{\chi}$ such that
$[c_{\chi},\i)\cap\left(\cup_{k=1}^{\i} A_{\chi, k}\right)$ is dense in $[c_{\chi},\i)$.
\end{theorem}
Let us note that Theorem \ref{den} shows in particular that $L(s, \chi)$ is not logarithmically completely  monotonic on any subinterval of $[c_{\chi},\i)$. A stronger assertion is as follows:

For any subinterval of $[c_{\chi},\infty)$, however small it may be, infinitely many derivatives $F^{(k)}(s,\chi)$  change sign in this subinterval.\\

Now consider any two points $s_1, s_2$ with $1<s_1<s_2$. In the case of the Riemann zeta function, if we compare the signs of the values of $\frac{d^k}{ds^k}\log\zeta(s)$ at $s_1$ and $s_2$ for all values of $k$, we see that they are always the same. Then a natural question arises - what can we say if we make the same comparison in the case of a Dirichlet $L$-function? We will see below that the answer is completely different (actually it is as different as it could be). We first define a function $\psi_{\chi}$ for a real primitive Dirichlet character modulo $q$ as follows:\\
 
Let $\mathcal{B}:=\{g:\mathbb{N}\to \{-1,0,1\}\}$. Define an equivalence relation $\sim$ on $\mathcal{B}$ by $g\sim h$ if and only if $g(n)=h(n)$ for all $n$ large enough. Let $\hat {\mathcal{B}}=\mathcal{B}/\sim$. By abuse of notation, we define $\psi_{\chi}:(1,\infty) \to \hat{\mathcal{B}}$ to be a function whose image is a sequence given by $\{\operatorname{sgn } (F^{(k)}(s,\chi))\}$, i.e.,
\begin{equation}
\psi_{\chi}(s)(k):=\operatorname{sgn } (F^{(k)}(s,\chi)). 
\end{equation}
With this definition, we answer the above question in the form of the following theorem.
\begin{theorem}\label{inj1}
Let $\chi$ be a real primitive character modulo $q$ and let $\psi_{\chi}$ be defined as above. Then there exists a constant $C_{\chi}$ with the following property:\\

\textup{(a)} The Riemann hypothesis for $L(s, \chi)$ implies that $\psi_{\chi}$ is injective on $[C_{\chi}, \infty)$.\\

\textup{(b)} Let $\psi_{\chi}$ be injective on $[C_{\chi},\infty)$. Then there exists an effectively computable constant $D_{\chi}$ such that if all the nontrivial zeros $\r$ of $L(s,\chi)$ up to the height $D_{\chi}$ lie on the critical line Re $s=1/2$, then the Riemann Hypothesis for $L(s,\chi)$ is true.
\end{theorem}
\section{Proof of theorem \ref{den}}

First we will compute $F^{(k)}(s,\chi)$ in terms of the zeros of $L(s,\chi)$. The logarithmic derivative of $L(s,\chi)$ satisfies \cite[page. 83]{dav} 
\begin{equation}\label{llp}
F'(s,\chi)=\frac{L'(s,\chi)}{L(s,\chi)}=-\frac{1}{2}\log \frac{q}{\pi}-\frac{1}{2}\frac{\G'(s/2+b/2)}{\G(s/2+b/2)}+B(\chi)+\sum_{\r}\left(\frac{1}{s-\r}+\frac{1}{\r}\right),
\end{equation}
where $B(\chi)$ is a constant depending on $\chi$, 
\begin{equation}\label{b}
b=\twopartdef{1}{\chi(-1)=-1}{0}{\chi(-1)=1},
\end{equation}
and $\r=\b+i\g$ are the non trivial zeros of $L(s,\chi)$. 
Since $B(\overline{\chi})=\overline{B(\chi)}$ and $\chi$ is real, $B(\chi)$ is given by 
\begin{equation*}
B(\chi)=-\sum_{\r}\frac{1}{\r}=-2\sum_{\g >0}\frac{\b}{\b^2+\g^2} < \infty,
\end{equation*}
see \cite[page. 83]{dav}. Note that $B(\chi)$ is negative. The Weierstrass infinite product for $\G(s)$ is \cite[p.~73]{dav}
\begin{equation}\label{G}
\G(s)=\frac{e^{-\g s}}{s}\prod_{n=1}^{\infty}(1+s/n)^{-1}e^{s/n},
\end{equation}
with $s=0,-1,-2,\dots$ being its simple poles. The functional equation for $\G(s)$ is
\begin{align}
\G(s+1)=s\G(s)\label{F1}\\
\end{align}
where as the duplication formula for $\G(s)$ is
\begin{align}
\G(s)\G(s+1/2)=2^{(1-2s)}\pi ^{1/2}\G(2s),\label{F2}
\end{align}
see \cite[p.~73]{dav}. The following can be easily derived from (\ref{F1}), (\ref{F2}) and the logarithmic derivative of (\ref{G}):
\begin{align}
&\frac{1}{2}\frac{\G'(s/2)}{\G(s/2)}=-\frac{\g}{2}-\frac{1}{s}-\sum_{n=1}^{\infty}\left(\frac{1}{s+2n}-\frac{1}{2n}\right),\label{GZE}\\
&\frac{1}{2}\frac{\G'(s/2+1/2)}{\G(s/2+1/2)}=-\log (2)-\frac{\g}{2}-\sum_{n=0}^{\infty}\left(\frac{1}{s+2n+1}-\frac{1}{2n+1}\right),\label{GZO}
\end{align}
From (\ref{llp}), (\ref{b}), (\ref{GZE}) and (\ref{GZO}), we have
\begin{equation}\label{llpf}
F'(s,\chi)=-\frac{1}{2}\log \frac{q}{\pi}+b \log 2+\frac{\g}{2}+B(\chi)+\frac{1-b}{s}+\sum_{\r\neq 0}\left(\frac{1}{s-\r}+\frac{1}{\r}\right),
\end{equation}
where $\r$ runs through all the zeros of $L(s,\chi)$. The successive differentiation of (\ref{llpf}) gives for $k\geq 2$,
\begin{align}\label{sd}
F^{(k)}(s,\chi)&=(-1)^{k-1} (k-1)!\left(\frac{1-b}{s^{k}}+\sum_{\substack{\r\neq 0\\L(\r, \chi)=0}}\frac{1}{(s-\r)^k}\right)\nonumber\\
&=(-1)^{k-1} (k-1)!\left(\sum_{L(\r, \chi)=0}\frac{1}{(s-\r)^k}\right).
\end{align}
\begin{center}
\begin{figure}[h!]
\begin{tikzpicture}[scale=3]
  \tkzDefPoint(0,0){O}
  \tkzDefPoint(1,0){1}
  \tkzDefPoint(.5,0){1/2}
	\tkzDefPoint(0,1.5){A1}
 \tkzDefPoint(1.5,0){AA}
	\tkzDefPoint(.5,1.5){B1}
  \tkzDefPoint(.5,0){B2}
  \tkzDefPoint(1,1.5){C1}
	\tkzDefPoint(1,0){C2}
	\tkzDefPoint(2,0){s_0}
	\tkzDefPoint(2.5,0){s_1}
	\tkzDefPoint(2.2,0){s_2}
	\tkzDefPoint(.8,.85){r_0}
	\tkzDefPoint(.9375,1.09){r_1}
	\tkzDrawArc[rotate,color=red](s_0,r_0)(30)
  \tkzDrawArc[rotate,color=red](s_0,r_0)(-30)
	\tkzDrawArc[rotate,color=red](s_1,r_0)(-30)
	\tkzDrawArc[rotate,color=red](s_1,r_0)(30)
	\tkzDrawArc[rotate, style=dashed,color=red](s_2,r_0)(-30)
	\tkzDrawArc[rotate, style=dashed,color=red](s_2,r_0)(30)	
	\tkzDrawLines(A1,O)
	\tkzDrawLines[style=dashed](B1,B2)
	\tkzDrawLines(C1,C2)
  \tkzDrawLines[add = 0 and 1](O,AA)
	\tkzDrawSegments[color=blue](s_0,r_0)
	\tkzDrawSegments[color=blue](s_1,r_0)
	\tkzDrawSegments[color=blue](s_1,r_1)
	\tkzDrawSegments[style=dashed,color=blue](s_2,r_0)
  \tkzDrawPoints(O,1,1/2)
	\tkzDrawPoints(s_0,s_1,s_2,r_0,r_1)
  \tkzLabelPoints[left](O,1)
	\tkzLabelPoints[below](s_0,s_1,s_2)
	\coordinate[label=above:$\r_1$] (r_1) at (.9375,1.09);
	\coordinate[label=below:$\r_0$] (r_0) at (.8,.85);
	\coordinate[label=above:$\s{=}1/2$] (B1) at (.5,1.5);
 \end{tikzpicture}
	\caption{Construction for identifying the unique $\r_0$ at which $l(s)$ is attained for $s\in(s_0-\e,s_0+\e)$.}
\end{figure}
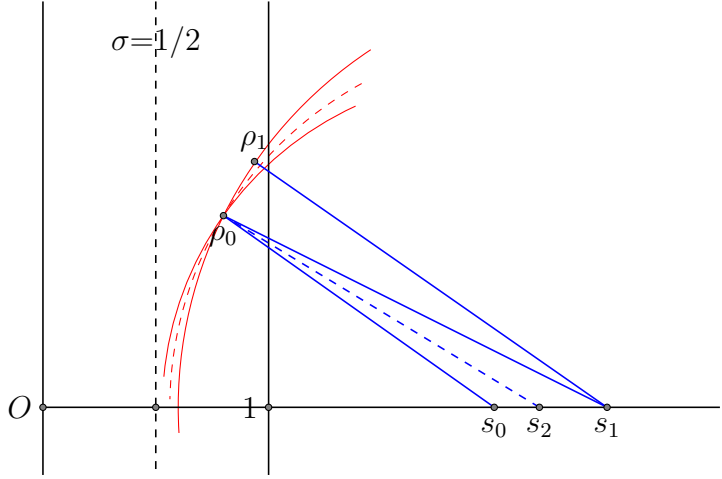
\end{center}
Let $s>1/2$ and define 
\begin{equation}\label{ells}
l(s):=\min\{|s-\r|:L(\r,\chi)=0\}.
\end{equation}
If the minimum $l(s)$ is attained for a non-trivial zero $\r$ of $L(s,\chi)$, then since the non-trivial zeros are symmetric with respect to the line $\sigma=1/2$, we have $\Re\r\geq 1/2$. Let $\tilde{\r}_0$ be the non-trivial zero of $L(s,\chi)$ with minimum but positive imaginary part, i.e., $\Im \tilde{\r}_0=\min\{\Im \r>0: L(\r,\chi)=0, \Re\r\geq 1/2 \}$. Write $\tilde{\r}_0=\tilde{\b}_0+i\tilde{\g}_0$. Then for  all $s>\tilde{\g}_0^2+1/4$, we have $s^2>(s-1/2)^2+\tilde{\g}_0^2\geq |s-\tilde{\r}_0|^2\geq (l(s))^2$. Define 
\begin{equation}\label{cchi}
c_{\chi}:=\inf\{c>1:s>c\Rightarrow |s|>l(s)\}.
\end{equation}
The constant $c_{\chi}$ is defined in this way since we want $l(s)$ to be attained at a non-trivial zero of $L(s,\chi)$, as this will allow us to separate the two terms of the series in (\ref{sd}) corresponding to this zero and its conjugate, which together will give a dominating term essential in the proof. Note that if $\tilde{\g}_0\leq\sqrt{3}/2$, $c_{\chi}=1$, otherwise $1 \leq c_{\chi}\leq \tilde{\g}_0^2+1/4$.

Next we show that for any $s\geq c_{\chi}$, there is an $s'\in (s-\e,s+\e)$, $\epsilon>0$, so that $l(s')$ is attained at a unique non-trivial zero $\r'$ of $L(s,\chi)$ with $\Im \r'>0$.

For any real number $s_0 > c_{\chi}$, consider the interval $(s_0-\e,s_0+\e) \subset [c_{\chi},\infty)$ for some $\e>0$.
Let 
\begin{equation}
A:=\{\r':\Im\r'\geq 0 \mbox{ and } |s_0-\r'|=l(s_0), L(\r',\chi)=0\},
\end{equation}
that is, $A$ is comprised of all non-trivial zeros on the circle with center $s_0$ and radius $l(s_0)$. Clearly $A$ is a finite set since $|A|\leq \mbox{N}(l(s_0),\chi)$, where $N(T,\chi)$ denotes the number of zeros of $L(s,\chi)$ up to height $T$. As shown in Figure 1, let $\r_0 \in A$ with $\Re\r_0=\max\{\Re \r: \r \in A\}$. Then for any $s\in (s_0,s_0+\e)$, $|s-\r_0|<|s-\r|$, for all $\r \in A$, $\r \neq \r_0$. Fix one such $s$, say $s_1$, so that $s_0<s_1<s_0+\e$. Now more than one zeros may lie on the circle with center $s_1$ and radius $|s_1-\r_0|$. If there aren't any (apart from $\r_0$), then we have constructed $s'(=s_1)$ that we sought. If there are more than one, we select the one among them, say $\r_1$, which has the minimum real part, i.e., $\Re \r_1=\min\{\Re \r:|s_1-\r_0|=|s_1-\r|, \r\neq\r_0, L(\r,\chi)=0\}$. Note that $\Im \r_1>\Im \r_0$, otherwise it will contradict the fact that the minimum $l(s_0)$ is attained at $\r_0$. 

For any $s\in (s_0,s_1)$, $|s-\r_0|<|s-\r_1|$. Now fix one such $s$, say $s_2\in (s_0,s_1)$, and find a $\r_2$ so that $\Re \r_2=\min\{\Re \r:|s_2-\r_0|=|s_2-\r|, \r\neq\r_0, L(\r,\chi)=0\}$. Since there are only finite may zeros in the rectangle $[0,1]\times[\Im \r_0,\Im \r_1]$, repeating the argument allows us to find an $s' \in \mathbb{R}$ and $s_0<s'<s_1<s_0+\e$, so that $\r_0$ is the only non-trivial zero of $L(s,\chi)$ with $\Im \r_0\geq 0$ and $|s'-\r_0|=\min\{|s'-\r|, \Im \r\geq 0\mbox{ and }L(\r,\chi)=0\}$, i.e., the circle with center $s'$ and radius $|s-\r_0|$ does not contain any zero other than $\r_0$ itself. Note that for any $s\in(s_0, s')$, $\r_0$ is the only zero at which $l(s)$ is attained.

Next, let $B=\{\r':\r'\neq\r_0 , |s_0-\r'|<|s_0-\r|\}$, where $\r$, $\r'$ are zeros of $L(s,\chi)$. Note that $B$ is also a finite set. Arguing in a similar way as above, we can find a $\tilde{\r}\in B$ and $s''\in (s_0, s_0+\e)$ so that for all $s\in (s_0, s'')$, $|s-\tilde{\r}|\leq|s-\r|$ for $\r\neq \r_0$.

Therefore we can find a closed interval $[c,d]\subset (s_0-\e,s_0+\e)$ so that for all $s \in [c,d]$, we have
\begin{align}
l(s)=|s-\r_0|=|s-\bar{\r_0}|<|s-\r|, \r\neq \r_0, \overline{\r_0}\label{ineq1}\\
|s-\tilde{\r}|=|s-\overline{\tilde{\r}}|\leq|s-\r|,  \r\neq \r_0, \overline{\r_0}, \tilde{\r}, \overline{\tilde{\r}}.\label{ineq2}
\end{align}

Now let $s-\r_0=r_s e^{i\t_s}$ for all $c\leq s \leq d$. Then from (\ref{sd}) and the fact that the zeros of $L(s,\chi)$ are symmetric with respect to the real axis, we have
{\allowdisplaybreaks
\begin{align}\label{sdl}
F^{(k)}(s,\chi)&=(-1)^{k-1} (k-1)!\left(\frac{1}{(s-\r_0)^k}+\frac{1}{(s-\bar{\r_0})^k}+\sum_{\r\neq\r_0,\bar{\r_0}}\frac{1}{(s-\r)^k}\right)\nonumber\\
&=(-1)^{k-1} (k-1)!\left(\frac{2}{r_s^k}\cos (k\t_s)+\sum_{\r\neq\r_0,\bar{\r_0}}\frac{1}{(s-\r)^k}\right)\nonumber\\
&=\frac{(-1)^{k-1} (k-1)!}{r_s^k}\left(2\cos (k\t_s)+f(s)\right),
\end{align}}
where $f(s):=r_s^k\sum_{\r\neq\r_0,\bar{\r_0}}\frac{1}{(s-\r)^k}$ and $k\geq 2$. Since the series $\sum_{\r\neq\r_0,\bar{\r_0}}\frac{1}{(s-\r)^k}$ converges absolutely for $k\geq 1$, $f(s)$ is a differentiable function for $s> 1$. Now,
{\allowdisplaybreaks
\begin{align}\label{est1}
|f(s)|\leq 2 \sum_{\substack{\r\neq\r_0,\bar{\r_0}\\\Im \r\geq 0}}\frac{r_s^k}{|s-\r|^k}
&=2\sum_{\substack{\r\neq\r_0,\bar{\r_0},\\\Im \r\geq 0}}\frac{|s-\r_0|^k}{|s-\r|^k}\nonumber\\
&= 2 |s-\r_0|^2\sum_{\substack{\r\neq\r_0,\bar{\r_0},\\\Im \r\geq 0}}\frac{1}{|s-\r|^2}\frac{|s-\r_0|^{k-2}}{|s-\r|^{k-2}}\nonumber\\
&\leq 2 |s-\r_0|^2\sum_{\substack{\r\neq\r_0,\bar{\r_0},\\\Im \r\geq 0}}\frac{1}{|s-\r|^2}\frac{|s-\r_0|^{k-2}}{|s-\tilde{\r}|^{k-2}}\nonumber\\
&\leq 2 |s-\r_0|^2\sum_{\substack{\r\neq\r_0,\bar{\r_0},\\\Im \r\geq 0}}\frac{1}{|s-\r|^2}\sup_{c\leq s\leq d}\left\{\frac{|s-\r_0|^{k-2}}{|s-\tilde{\r}|^{k-2}}\right\},
\end{align}}%
where in the penultimate step we use (\ref{ineq2}). Let $h(s):=\frac{|s-\r_0|}{|s-\tilde{\r}|}$. Then $h(s)$ is a continuous function on $[c,d]$ and hence attains its supremum on $[c,d]$. Thus there exists an $x\in [c,d]$ such that 
\begin{equation}\label{eta}
\eta:=\sup_{c\leq s\leq d}\left\{\frac{|s-\r_0|}{|s-\tilde{\r}|}\right\}=\frac{|x-\r_0|}{|x-\tilde{\r}|}.
\end{equation}
Therefore by (\ref{ineq1}), $\eta<1$. Combining (\ref{est1}) and (\ref{eta}), we have
\begin{align}\label{est2}
|f(s)|&\leq 2\eta^{k-2} |s-\r_0|^2\sum_{\substack{\r\neq\r_0,\bar{\r_0}\\\Im \r\geq 0}}\frac{1}{|s-\r|^2}\leq 2\eta^{k-2} |d-\r_0|^2\sum_{\substack{\r\neq\r_0,\bar{\r_0}\\\Im \r\geq 0}}\frac{1}{|c-\r|^2}\leq C_{c,d,\chi}\eta^{k-2}.
\end{align}%
Note that the constant term depends only on $c,d$ and $\chi$. Hence for sufficiently large $k$, we have $|f(s)|<1$. Let $c-\r_0=r_c e^{i \t_c}$ and $d-\r_0=r_d e^{i \t_d}$. Then $\t_c>\t_d$. For $k$ large enough, we can write $2\pi<k(\t_c-\t_d)$. Since for $s\in [c,d]$, we have $\t_d\leq\t_s\leq\t_c$, for a sufficiently large $k$, $\cos(k\t_s)$ attends all the values of the interval $[-1,1]$. So from (\ref{est1}) and (\ref{est2}) we conclude that for each large enough $k$ there will be an $s$ in $[c,d]\subset (s_0-\e,s_0+\e)$ so that $F^{(k)}(s,\chi)=0$. This shows that $\cup_{k=1}^{\i} A_{\chi, k}$ has a non-empty intersection with $(s_0-\e,s_0+\e)$ for any $s_0>c_{\chi}$. This completes the proof of the theorem.\\

\textbf{Remark:} Let $\chi$ be a real nonprincipal Dirichlet character. 
If $L(s,\chi)$ has a Siegel zero, call it $\b$, and if every zero of $L(s,\chi)$ has real part $\leq \b$, then for any $s>1$, $(\ref{sd})$ implies
\begin{align}\label{sdi}
F^{(k)}(s,\chi)=\frac{(-1)^{k-1} (k-1)!}{(s-\b)^k}\left(1+\sum_{\r\neq \b \atop L(\r, \chi)=0}\left(\frac{s-\b}{s-\r}\right)^k\right).
\end{align}
Arguing as in the proof of Theorem \ref{den}, we see that there exists an integer $M$ such that for all $k\geq M$, the series in (\ref{sdi}) is less than $1$. This means that for those $k$, $F^{(k)}(s,\chi)$ maintains the same sign for all $s>1$. This is why we include the condition that $L(s,\chi)\neq 0$ for $0<s<1$ in the hypotheses of Theorem \ref{den} . 
\section{Proof of theorem \ref{inj1}}

Assume that the Riemann hypothesis holds for $L(s,\chi)$. Let $\g_0:=\Im \r_0=\min\{\Im \r\geq 0: L(\r,\chi)=0\}$, where $\r_0,\r$ are non-trivial zeros of $L(s,\chi)$. Then $\r_0=1/2+i\g_0$. We show that the function $\psi_{\chi}$ is injective on $[C_{\chi},\infty)$, where the constant $C_{\chi}$ will be determined later. 

Let $s>c_{\chi}$, where $c_{\chi}$ is defined in (\ref{cchi}). Then $l(s)<|s|$ and $l(s)=|s-\r_0|<|s-\r|$ for $\r \neq \r_0, \bar{\r_0}$. Let $s-\r_0=r_s e^{i \t_s}$. From (\ref{sdl}), we have for $k\geq 2$,
\begin{align}\label{err}
|f(s)|\leq\sum_{\r\neq\r_0,\bar{\r_0}}\frac{r_s^k}{|s-\r|^k}&=|s-\r_0|^2\sum_{\r\neq\r_0,\bar{\r_0}}\frac{1}{|s-\r|^2}.\frac{|s-\r_0|^{k-2}}{|s-\r|^{k-2}}\nonumber\\
&\leq |s-\r_0|^2\sum_{\r\neq\r_0,\bar{\r_0}}\frac{1}{|s-\r|^2}.\sup_{\r}\left\{\frac{|s-\r_0|^{k-2}}{|s-\r|^{k-2}}\right\}\nonumber\\
&=|s-\r_0|^2 \eta_s^{k-2}\sum_{\r\neq\r_0,\bar{\r_0}}\frac{1}{|s-\r|^2}\nonumber\\
&=O_{s,\chi}(\eta_s^{k-2}).
\end{align}
Here in the penultimate step, 
\begin{equation*}
\eta_s=\sup_{\r}\left\{\frac{|s-\r_0|}{|s-\r|}\right\}\leq\frac{|s-\r_0|}{|s-\tilde{\r}|}< 1,
\end{equation*}
and $\Im \r_0<\Im\tilde{\r}\leq\Im \r$, resulting from (\ref{ineq1}) and (\ref{ineq2}). Combining (\ref{sdl}) and (\ref{err}), we obtain
\begin{equation}\label{injj1}
F^{(k)}(s,\chi)=\frac{2 (-1)^{k-1} (k-1)!}{r_s^k}\left(\cos (k\t_s)+f(s)\right),
\end{equation}
where $f(s)=O_{s,\chi}(\eta_s^{k-2})$.

Next, we show that there are infinitely many $k$ for which $\cos (k\t_s)$, which we view as the main term, dominate the error term. Since $\eta_s<1$, for a fixed $s>1$, we can bound the error term in $(-\e,\e)$ for all sufficiently large $k$ and for all $0<\e<1$. Write $\cos (k \t_s)=\cos \left(\pi \frac{k \t_s}{ \pi}\right)=\cos \left(2 \pi \frac{k \t_s}{2 \pi}\right)$ and consider the cases when $\frac{\t_s}{\pi}$ is rational and $\frac{\t_s}{2 \pi}$ is irrational. 

If $\frac{\t_s}{\pi}$ is a rational number, there are infinitely many $k\in \mathbb{N}$ so that $\frac{k\t_s}{2 \pi}$ is an even integer and hence $\cos (k \t_s)=1$.

If $\frac{\t_s}{\pi}$ is a rational number with odd numerator, then there are infinitely many $k\in \mathbb{N}$, namely the odd multiples of the denominator, so that $\frac{k\t_s}{2 \pi}$ is an odd integer and hence $\cos (k \t_s)=-1$.

Let $\frac{\t_s}{\pi}=\frac{2m}{n}$ be a rational number with even numerator and odd denominator. Since $(2m,n)=1$, there exists an integer $l\in [1,n]$ such that $2ml\equiv 1 (\operatorname{mod} n)$. For all $k\equiv l (\operatorname{mod} n) ,2mk\equiv 1 (\operatorname{mod} n)$. Therefore for all $k\equiv l \operatorname{mod} n $, since $2mk$ is even, we have $2mk=(2p+1)n+1$.
Hence there are infinitely many integers $k$ for which $\cos(k\t_s)=\cos\left(\pi \left(2p+1+\frac{1}{n}\right)\right)=-\cos\left(\frac{\pi}{n}\right)$.

If $\frac{\t_s}{2 \pi}$ is irrational, then we know from \cite{wey} that the sequence $\left\{\left\{\frac{k\t_s}{2 \pi}\right\}\right\}$ is dense in $[0,1]$, where $\{x\}$ denotes the fractional part of $x$. (Actually, Kronecker's approximation theorem is sufficient to prove the denseness.) Hence there are infinitely many $k\in \mathbb{N}$ with $\left\{\frac{k\t_s}{2 \pi}\right\}$ close to $1$ and hence $\cos (k \t_s)>1-\e$ for any given $\e>0$. Likewise, there are infinitely many $k\in \mathbb{N}$ with $\left\{\frac{k\t_s}{2 \pi}\right\}$ close to $\frac{1}{2}$ and hence $\cos (k \t_s)<-1+\e$.

Fix $s_1$ and $s_2$ such that $c_{\chi}<s_1<s_2$. Then $l(s_1)=|s_1-\r_0|$ and $l(s_2)=|s_2-\r_0|$. Let $\t_1$ and $\t_2$ be such that $s_1-\r_0=r_1 e^{i \t_1}$ and $s_2-\r_0=r_2 e^{i \t_2}$. Note that $0<\t_2<\t_1<\pi/2$. From (\ref{injj1}), we have
\begin{align}
F^{(k)}(s_1,\chi)&=\frac{2 (-1)^{k-1} (k-1)!}{r_1^k}\left(\cos (k\t_1)+f(s_1)\right),\label{in1}\\
F^{(k)}(s_2,\chi)&=\frac{2 (-1)^{k-1} (k-1)!}{r_2^k}\left(\cos (k\t_2)+f(s_2)\right)\label{in2},
\end{align}
where $f(s_1)=O_{s_1,\chi}(\eta_{s_1}^{k-2})$ and $f(s_2)=O_{s_2,\chi}(\eta_{s_2}^{k-2})$.
Write $\t_1=\t_2+(\t_1-\t_2)$.

We show that there exist infinitely many integers $k$ such that the terminal rays of $k\t_1$ and $k\t_2$ stay away from the $y$-axis, that $\text{sgn}\left(\cos k\t_1\right)=-\text{sgn}\left(\cos k\t_2\right)\neq 0$, and that $\cos(k\t_1)$ and $\cos(k\t_2)$ dominate $f(s_1)$ and $f(s_2)$ in (\ref{in1}) and (\ref{in2}) respectively. We first determine the signs.

\underline{Case 1:} If $\frac{\t_1-\t_2}{\pi}$ is rational with odd numerator then as we saw before, there are infinitely many positive integers $k$ so that $k\frac{(\t_1-\t_2)}{\pi}$ is an odd integer and hence for those $k \in \mathbb{N}$, $\cos(k \t_1)=\cos(k\t_2+\pi)=-\cos(k \t_2)$.

\underline{Case 2:} If $\frac{\t_1-\t_2}{\pi}$ is rational with even numerator and odd denominator $n$, there are infinitely many positive integers $k$ so that $k\frac{(\t_1-\t_2)}{\pi}=2p+1+1/n$ for some $p \in \mathbb{N}$ and so $\cos(k \t_1)=\cos(k\t_2+\pi+\pi/n )=-\cos(k\t_2+\pi/n)$. 

\underline{Case 3:} If $\frac{(\t_1-\t_2)}{2\pi}$ is irrational, there are infinitely many positive integers $k$ so that $\left\{k\frac{(\t_1-\t_2)}{2\pi}\right\}\in [1/2,1/2+\e/2\pi)$, for any given $\e>0$. So for any $\delta$ such that $0<\d<\e$, we have $\cos(k\t_1)=\cos(k \t_2+\pi+\d)=-\cos(k\t_2+\d)$. We can choose $\e$ as small as we want and hence $0<\d<\e<\pi/n$. 

We first show that in Case 2, we have the terminal rays of the angles sufficiently away from the $y$-axis, with $\cos k\t_1$ and $\cos k\t_2$ dominating their corresponding terms $f(s_1)$ and $f(s_2)$. To that end, choose a constant $b_{\chi}>1/2$ such that $\tan\left(\frac{\pi}{100}\right)=\frac{\g_0}{b_{\chi}-\frac{1}{2}}$, say. If $s-\r_0=r_se^{i\t_s}$ and $s>b_{\chi}$, then $0<\t_s<\pi/100$. So if we take $b_{\chi}<s_1<s_2$, then $0<\t_2<\t_1<\pi/100$. Since $\eta_{s_1}, \eta_{s_2}<1$ there exists an integer $K$ such that $|f(s_1)|, |f(s_2)|<\t_2/4$ for all $k>K$. As we saw before, for infinitely many integers $k>K+2$, we have $k\t_1=k\t_2+\pi+\pi/n$, where $n$ depends on $\t_1$ and $\t_2$. We first note that all angles below are considered mod $2\pi$. If $k\t_2 \in (\pi/2+\t_2, \pi)$ then $k \t_1\in (-\pi/2+\t_2, \pi/2-\t_2)$. Thus $\cos(k\t_1)\cos(k\t_2)<0$. Also $|\cos(k\t_2)|>|\sin(\t_2)|\geq \t_2/2>|f(s_2)|$ and $|\cos(k\t_1)|=|\cos(k\t_2+\pi/n)|>|\sin(\t_2)|\geq \t_2/2>|f(s_1)|$.

Similarly we see that $|\cos(k\t_1)|>|f(s_1)|$ and $|\cos(k\t_2)|>|f(s_2)|$ when $k\t_2 \in (-\pi/2+\t_2,0) $. If $k\t_2 \in (0,\pi/2-\t_2)$ and $k \t_1\in (-\pi, -\pi/2-\t_2)$ in this case also $|\cos(k\t_2)|>|\sin(\t_2)|\geq \t_2/2>|f(s_2)|$ and $|\cos(k\t_1)|=|\cos(k\t_2+\pi/n)|>|\sin(\t_2)|\geq \t_2/2>|f(s_1)|$. Now let $k\t_2 \in (0,\pi/2+\t_2)$ and  $k \t_1\in (-\pi/2-\t_2, 0)$. Then since $\pi/n <\t_1<\pi/100$, it is easy to check that $(k-2)\t_2 \in (0,\pi/2-\t_2)$ and $(k-2)\t_1=k\t_2+\pi+\pi/n-2\t_1 \in (-\pi, -\pi/2-\t_2)$. Hence $|\cos(k\t_1)|>|f(s_1)|$ and $|\cos(k\t_2)|>|f(s_2)|$. Similarly we have the same conclusion if $k\t_2 \in (-\pi,-\pi/2+\t_2)$ and $k \t_1\in (\pi/2-\t_2, \pi)$.

Note that since $k\t_2+\pi+\pi/n>k\t_2+\pi+\d$, for the values of $\t_1$ and $\t_2$ in Case 3 as well, one can similarly prove that $|\cos(k\t_2)|>|f(s_2)|$ and $|\cos(k\t_1)|>|f(s_1)|$. So is the case with the values of $\t_1$ and $\t_2$ in Case 1. 

Let 
\begin{equation}\label{capcchi}
C_{\chi}=\max \{c_{\chi},b_{\chi}\}.
\end{equation}
Then for any given real numbers $s_1$ and $s_2$ such that $C_{\chi}<s_1<s_2$, we have shown that there exist infinitely many integers $k$ such that $\cos(k\t_1)$ and $\cos(k\t_2)$ have opposite signs and $|\cos(k\t_1)|>|f(s_1)|$ and $\cos(k\t_2)>f(s_2)$. This implies that $F^{(k)}(s_1,\chi)$ and $F^{(k)}(s_2,\chi)$ have opposite signs and that in turn proves that the function $\psi_{\chi}$ is injective in $[C_{\chi},\infty)$.\\

We now prove part (b) of Theorem \ref{inj1}.
\begin{center}
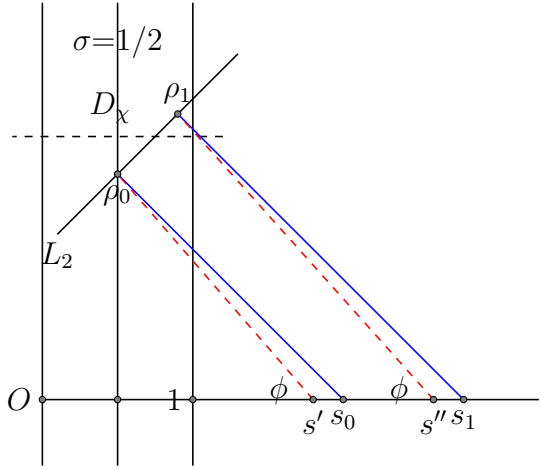
\begin{figure}[h!]
\begin{tikzpicture}[scale=2]
  \tkzDefPoint(0,0){O}
  \tkzDefPoint(1,0){1}
  \tkzDefPoint(.5,0){1/2}
	\tkzDefPoint(0,2.2){A1}
 \tkzDefPoint(1.5,0){AA}
	\tkzDefPoint(.5,2.2){B1}
  \tkzDefPoint(.5,0){B2}
  \tkzDefPoint(1,2.2){C1}
	\tkzDefPoint(1,0){C2}
	\tkzDefPoint(2,0){s_0}
	\tkzDefPoint(1,1.75){AB}
	\tkzDefPoint(0,1.75){BA}
	\tkzDefPoint(2.8,0){s_1}
	\tkzDefPoint(.5,1.5){r_0}
	\tkzDefPoint(.9,1.9){r_1}
		\tkzDefPoint(1.8,0){s'}
	\tkzDefPoint(2.6,0){s''}
	\tkzDrawLines(A1,O)
	\tkzDrawLines(B1,B2)
	\tkzDrawLines[style=dashed](AB,BA)
	\tkzDrawLines(C1,C2)
	\tkzDrawLines[add = 1 and 1](r_0,r_1)
  \tkzDrawLines[add = 0 and 1.2](O,AA)
	\tkzDrawSegments[color=blue](s_0,r_0)
	\tkzDrawSegments[color=blue](s_1,r_1)
	\tkzDrawSegments[color=red,style=dashed](s',r_0)
	\tkzDrawSegments[color=red,style=dashed](s'',r_1)
  \tkzDrawPoints(O,1,1/2)
	\tkzDrawPoints(s_0,s_1,r_0,r_1,s',s'')
  \tkzLabelPoints[left](O,1)
	\tkzLabelPoints[below](s_0,s_1,s',s'')
	\coordinate[label=above:$\r_1$] (r_1) at (.9,1.9);
	\coordinate[label=below:$\r_0$] (r_0) at (.5,1.5);
	\coordinate[label=above:$\s{=}1/2$] (B1) at (.5,2.2);
	\coordinate[label=above:$D_{\chi}$] (B1) at (.45,1.77);
	\coordinate[label=above:$\phi$] (s') at (1.57,-.1);
	\coordinate[label=above:$\phi$] (s'') at (2.37,-.1);
	\coordinate[label=above:$L_2$] (s) at (.1,.8);
 \end{tikzpicture}
	\caption{Constructing the angle $\phi=2\pi(a+b\sqrt{2})$.}
\end{figure}
\end{center} 
Let $\r_0$ be  the lowest zero of $L(s,\chi)$ above the real axis (so $\r_0$ is not a real number).
Let $L_1$ be the line passing through $\r_0$ and perpendicular to the line which passes through $\r_0$ and $C_{\chi}$, where $C_{\chi}$ is defined in (\ref{capcchi}). Let $(1, D_{\chi})$ be the point of intersection of the lines $\sigma=1$ and $L_1$. We first show that if there is only one zero $\r_1$ with $\Im \r_1\geq D_{\chi}$ off the critical line $\sigma=1/2$, then this contradicts the injectivity of $\psi_{\chi}$ on $[C_{\chi},\infty)$. 

Without loss of generality, let $\Re \r_1>1/2$. As shown in Figure 2, let $L_2$ be the line passing through $\r_0$ and $\r_1$. Let $s_0$ and $s_1$ be the points of intersection of the real axis with the lines perpendicular to $L_2$ and passing through $\r_0$ and $\r_1$ respectively. Clearly $s_1>s_0>C_{\chi}$. Note that by our construction, $l(s_0)=|s_0-\r_0|$ and $l(s_1)=|s_1-\r_1|$, where $l(s)$ is defined in (\ref{ells}), and there exists a $\t$ such that $(s_0-\r_0)=r_{s_0}e^{i\t}$ and $(s_1-\r_1)=r_{s_1}e^{i\t}$.
From the proof of the Theorem \ref{den}, we know that there exists an $\e>0$ so that $l(s)=|s-\r_0|$ for all $s\in (s_0-\e,s_0+\e)$ and $l(s)=|s-\r_1|$ for all $s\in (s_1-\e,s_1+\e)$. Without loss of generality, we can assume that $s_0+\e<s_1-\e$. Therefore, there exists a $\d>0$ such that $\t_s\in (\t-\d,\t+\d)$, where $s-\r_0=r_se^{i\t_s}$ and $l(s)=|s-\r_0|$ for all $s\in (s_0-\e,s_0+\e)$, and such that $\t_s\in (\t-\d,\t+\d)$, where $s-\r_1=r_se^{i\t_s}$ and $l(s)=|s-\r_1|$ for all $s\in (s_1-\e,s_1+\e)$.

Since the sequence $\{\{n\sqrt{2}\}\}$ is dense in $[0,1)$, and $\{n\sqrt{2}\}=n\sqrt{2}-\lfloor n\sqrt{2}\rfloor$, there exists an integer $a$ and an integer $b\neq 0$ such that $a+b\sqrt{2}\in(\frac{\t-\d}{2\pi},\frac{\t+\d}{2\pi})$. Let $\phi=2\pi(a+b\sqrt{2})$, $s'\in (s_0-\e,s_0+\e)$ and $s''\in (s_1-\e,s_1+\e)$ be such that  $s'-\r_0=r_{s'}e^{i\phi}$ and $s''-\r_1=r_{s''}e^{i\phi}$. Therefore,
\begin{align}
F^{(k)}(s',\chi)&=\frac{2 (-1)^{k-1} (k-1)!}{r_{s'}^k}\left(\cos (k\phi)+f(s')\right)\label{in11}\\
F^{(k)}(s'',\chi)&=\frac{2 (-1)^{k-1} (k-1)!}{r_{s''}^k}\left(\cos (k\phi)+f(s'')\right)\label{in22},
\end{align}
where $|f(s')|=O(\eta_{s'}^{k-2})$ and $|f(s'')|=O(\eta_{s''}^{k-2})$. Let $\eta=\min\{\eta_{s'},\eta_{s''}\}$. Then $|f(s')|,|f(s'')|\leq C_{s',s''}\eta^{k-2}$ for some constant $C_{s',s''}$.

We next show that there exist positive constants $C_{a,b}$ and $K_{a,b}$ so that 
\begin{equation}
|4k(a+b\sqrt{2})+r|>\frac{C_{a,b}}{k},
\end{equation}
for any integers $r$ and $k$, with $k>K_{a,b}$. Let $|4k(a+b\sqrt{2})+r|\leq 1$. Then,
\begin{equation}
|4k(a-b\sqrt{2})+r|\leq |4k(a+b\sqrt{2})+r|+8k|b|\sqrt{2}\leq 1+8k|b|\sqrt{2}< \frac{k}{C_{a,b}}.
\end{equation}
Therefore for $k\geq2$,
\begin{equation}
|4k(a+b\sqrt{2})+r|\frac{k}{C_{a,b}}> |4k(a-b\sqrt{2})+r||4k(a+b\sqrt{2})+r|=|(4ka+r)^2-2(4kb)^2|\geq 1,
\end{equation}
since $b\neq 0$. If $|4k(a+b\sqrt{2})+r|\geq 1$, then of course, there exists a $K_{a,b}$, such that for $k>K_{a,b}$, we have $|4k(a+b\sqrt{2})+r|>\frac{C_{a,b}}{k}$. 
Hence in conclusion, for a large positive integer $N$ and for all $k>N$, if we choose $m$ so that $|4k(a+b\sqrt{2})\pm 1\pm 4m|<1$, we have
\begin{align}
|\cos k\phi|=\left|\sin\frac{\pi}{2}(4k(a+b\sqrt{2})\pm 1\pm 4m)\right|&\geq \sin\left(\frac{\pi C_{a,b}}{2k}\right)\nonumber\\
&\geq\frac{\pi C_{a,b}}{4k}\nonumber\\
&>C_{s',s''}\eta^{k-2}.
\end{align}
Therefore for the above mentioned $s'$ and $s''$ such that $s'\neq s''$, and for all $k>N$, $F^{(k)}(s',\chi)$ and $F^{(k)}(s'',\chi)$ have the same sign. This contradicts the injectivity of $\psi_{\chi}$ on $[C_{\chi},\infty)$. Now if there is more than one zero $\r$ with $\Im{\r}\geq D_{\chi}$ off the critical line, then we can choose the zero $\r_1$ with the following properties:

\textit{i}) The angle between the positive $x$-axis and the line $L$ passing through the zeros $\r_0$ and $\r_1$ is smaller than the angle between the positive $x$-axis and the line passing through the zeros $\r_0$ and $\r\neq\r_1$ and,

\textit{ii}) $\Im{\r_1}=\min\{\Im\r\geq D_{\chi}:\r\text{ lies on the line } L\}$.

 Then we can proceed similarly as above and again get a contradiction. Hence, all the zeros above the line $t=D_{\chi}$ lie on the critical line $\sigma=1/2$. This completes the proof.
\begin{center}
\textbf{Acknowledgements}
\end{center}
The first author is funded in part by the grant NSF-DMS 1112656 of Professor Victor H. Moll of Tulane University and sincerely thanks him for the support. 

\end{document}